\documentclass[final,1p,times]{elsarticle}

\usepackage{amssymb}
 \usepackage{amsthm}
\usepackage{amscd}
\usepackage{amsmath}
\usepackage{amsfonts}
\usepackage{amssymb}
\usepackage{graphicx}

\numberwithin{equation}{section}

\newcommand{\ra}{\rightarrow}
\newcommand{\p}{\partial}
\newcommand{\f}{\frac}

\newcommand{\be}{\begin{equation}}
\renewcommand{\ra}{\rightarrow}
\newcommand{\ee}{\end{equation}}
\newcommand{\bea}{\begin{eqnarray}}
\newcommand{\eea}{\end{eqnarray}}
\newcommand{\bna}{\begin{eqnarray*}}
\newcommand{\ena}{\end{eqnarray*}}

\renewcommand{\le}{\left}
\newcommand{\ri}{\right}

\journal{***}

\begin{document}

\begin{frontmatter}

\title{Existence of positive solutions to quasi-linear elliptic equations with exponential growth
in the whole Euclidean space}

\author{Yunyan Yang}
 \ead{yunyanyang@ruc.edu.cn}
\address{ Department of Mathematics,
Renmin University of China, Beijing 100872, P. R. China}

\begin{abstract}
In this paper a quasi-linear elliptic equation in the whole
Euclidean space is considered. The nonlinearity of the equation is
assumed to have exponential growth or have critical growth in view
of Trudinger-Moser type inequality. Under some assumptions on the
potential and the nonlinearity, it is proved that there is a
nontrivial positive weak solution to this equation. Also it is shown
that there are two distinct positive weak solutions to a
perturbation of the equation. The method of proving these results is
combining Trudinger-Moser type inequality, Mountain-pass theorem and
Ekeland's variational principle.

\end{abstract}

\begin{keyword}
Trudinger-Moser inequality\sep singular Trudinger-Moser inequality\sep
$N$-Laplace equation\sep exponential growth

\MSC 46E30\sep 46E35\sep 35J20\sep 35J60\sep 35B33

\end{keyword}

\end{frontmatter}

%\tableofcontents
\section{Introduction}
Let $\Omega\subset\mathbb{R}^N$ be a bounded smooth domain. There are fruitful
results on the following problem
\be\label{sob}\le\{\begin{array}{lll}-\Delta_pu=f(x,u)\quad{\rm in}\quad\Omega\\[1.5ex]
u\in W_0^{1,p}(\Omega),\end{array}\ri.\ee where $\Delta_pu={\rm
div}(|\nabla u|^{p-2}\nabla u)$. When $p=2$ and $|f(x,u)|\leq
c(|u|+|u|^{q-1})$, $1<q\leq 2^*=2N/(N-2)$, $N\geq 3$. Among pioneer
works we mention Br\'ezis \cite{Brezis}, Br\'ezis-Nirenberg
\cite{B-N}, Bartsh-Willem \cite{B-W} and Capozzi-Fortunato-Palmieri
\cite{CFP}. For $p\leq N$ and $p^2\leq N$, Garcia-Alonso \cite{G-A}
generalized Br\'ezis-Nirenberg's existence and nonexistence results
to $p$-Laplace equation. When $\Omega=\mathbb{R}^N$ and $p=2$, one
may consider the semilinear Schr\"odinger equation instead of
(\ref{sob}):
\be\label{sob-1}\le\{\begin{array}{lll}-\Delta u+V(x)u=f(x,u)\quad{\rm in}\quad\mathbb{R}^N\\[1.5ex]
u\in W^{1,N}(\mathbb{R}^N),\end{array}\ri.\ee where again
$|f(x,u)|\leq c(|u|+|u|^{q-1})$, $1<q\leq 2^*=2N/(N-2)$. Many papers
are devoted to (\ref{sob-1}), we refer the reader to
Kryszewski-Szulkin \cite{K-S}, Alama-Li \cite{A-L}, Ding-Ni
\cite{D-N} and Jeanjean \cite{Jea}.
 Sobolev embedding theorem and the critical point theory,
particularly moumtain-pass theorem would play an important role in
studying problems (\ref{sob}) and (\ref{sob-1}) since both of them
have variational structure. When $p=N$ and $f(x,u)$ behaves like
$e^{\alpha |u|^{{N}/{(N-1)}}}$ as $|u|\ra\infty$, problem (1.1) was
studied by Adimurthi \cite{Adi1}, Adimurthi-Yadava\cite{Adi3}, Ruf
et al \cite{ddR,dMR}, J. M. do \'O \cite{doo2}, Panda \cite{Panda2}
and the references therein. To the author's knowledge, all theses
results are based on Trudinger-Moser inequality
\cite{Moser,Pohozave,Trudinger} and critical point
theory.\\

In this paper we consider the existence of positive solutions of the
 quasi-linear equation \be\label{prob-0}
 -\Delta_Nu+V(x)|u|^{N-2}u=\f{f(x,u)}{|x|^\beta},\,\,\,x\in
  \mathbb{R}^N\,\,(N\geq 2),
\ee where $\Delta_Nu={\rm div}(|\nabla u|^{N-2}\nabla u)$,
$V:\mathbb{R}^N\ra \mathbb{R}$ is a continuous function, $f(x,s)$ is
continuous in $\mathbb{R}^N\times\mathbb{R}$ and behaves like
$e^{\alpha s^{{N}/{(N-1)}}}$ as $s\ra +\infty$, and $0\leq\beta<N$.
Problem (\ref{prob-0}) can be compared with (\ref{sob-1}) in this
way: Sobolev embedding theorem can be applied to (\ref{sob-1}),
while Trudinger-Moser type embedding theorem can be applied to
(\ref{prob-0}). When $\beta=0$, problem (\ref{prob-0}) was studied
by D. Cao \cite{Cao} in the case $N=2$, by Panda \cite{Panda}, J. M.
do \'O \cite{doo1} and Alves-Figueiredo \cite{Alves} in general
dimensional case. When $0<\beta<N$, problem (\ref{prob-0}) is
closely
related to a singular Trudinger-Moser type inequality, namely\\

\noindent{\bf Theorem A (\cite{Adi-Yang}).} {\it For all $\alpha>0$,
$0\leq\beta<N$, and $u\in W^{1,N}(\mathbb{R}^N)$ $(N\geq 2)$, there
holds
 \be\label{ss}\int_{\mathbb{R}^N}\f{e^{\alpha|u|^{N/(N-1)}}-\sum_{k=0}^{N-2}
 \f{\alpha^k|u|\,^{kN/(N-1)}}{k!}}{|x|^\beta}dx<\infty.\ee
 Furthermore, we have for all $\alpha\leq
 \le(1-\f{\beta}{N}\ri)\alpha_N$ and $\tau>0$,
 \be\label{criticalineq}\sup_{\int_{\mathbb{R}^N}(|\nabla
  u|^N+\tau|u|^N)dx\leq 1}
 \int_{\mathbb{R}^N}\f{e^{\alpha|u|^{N/(N-1)}}-\sum_{k=0}^{N-2}
 \f{\alpha^k|u|\,^{kN/(N-1)}}{k!}}{|x|^\beta}dx<\infty.\ee
   This inequality is sharp : for any
$\alpha>\le(1-\f{\beta}{N}\ri)\alpha_N$, the supremum is
  infinity.}\\

  This theorem extends a result of Adimurthi-Sandeep \cite{A-S} on a bounded smooth domain. When $\beta=0$ and $\tau=1$,
  (\ref{criticalineq}) was proved by B. Ruf in the case $N=2$ via symmetrization method and by
  Li-Ruf \cite{Li-Ruf} in general dimensional case via the method of
  blow-up analysis. When $\beta=0$ and $\alpha<\alpha_N$, (\ref{criticalineq}) was first proved by Cao \cite{Cao} in the case
  $N=2$, and then by Panda \cite{Panda}, J. M. do \'O \cite{doo1} in general dimensional case. A similar but different
  type inequality was obtained by Adachi-Tanaka \cite{AT}.\\

 We assume the following two conditions on the potential $V(x)$: \\

  \noindent$(V_1)$ $V(x)\geq V_0>0$ in $\mathbb{R}^N$ for some $V_0>0$;\\[1.5ex]
$(V_2)$ The function $\f{1}{V(x)}$ belongs to
$L^{\f{1}{N-1}}(\mathbb{R}^N)$.\\

 As for the nonlinearity $f(x,s)$ we suppose the following:\\

\noindent $(H_1)$ There exist constants $\alpha_0$, $b_1$, $b_2>0$
such that for all $(x,s)\in \mathbb{R}^N\times\mathbb{R}^+$,
$$|f(x,s)|\leq
b_1
s^{N-1}+b_2\le\{e^{\alpha_0|s|^{N/(N-1)}}-S_{N-2}(\alpha_0,s)\ri\};$$

\noindent $(H_2)$ There exists $\mu>N$ such that for all $x\in
\mathbb{R}^N$ and $s>0$,
$$0<\mu F(x,s)\equiv \mu\int_0^sf(x,t)dt\leq sf(x,s);$$

\noindent $(H_3)$ There exist constant $R_0$, $M_0>0$ such that for
all $x\in\mathbb{R}^N$ and $s\geq R_0$,
$$F(x,s)\leq M_0f(x,s).$$

% Since we are interested in positive solutions, as J. M. do \'O etc. did in \cite{doo1,doo3},
% we may assume $f(x,s)=0$ for all $(x,s)\in \mathbb{R}^N\times(-\infty,0]$.
% Moreover we assume the following growth condition on the
% nonlinearity $f(x,s)$:\\

\noindent Define a function space
\be\label{space}E=\le\{u\in W^{1,N}(\mathbb{R}^N):\int_{\mathbb{R}^N}V(x)|u|^Ndx<\infty\ri\}.\ee
We say that $u\in E$ is a weak solution of problem (\ref{prob-0}) if
for all $\varphi\in E$ we have
$$\int_{\mathbb{R}^N}\le(|\nabla u|^{N-2}\nabla u\nabla \varphi+V(x)|u|^{N-2}u\varphi\ri)dx
=\int_{\mathbb{R}^N}\f{f(x,u)}{|x|^\beta}\varphi dx.$$ The
assumption $(V_1)$ implies that $E$ is a
reflexive Banach space when equipped with the norm
\be\label{E-norm}\|u\|_E\equiv \le\{\int_{\mathbb{R}^N}\le(|\nabla
u|^N+V(x)|u|^N\ri)dx\ri\}^{\f{1}{N}}\ee and for any $q\geq N$, the
embedding
$$E\hookrightarrow W^{1,N}(\mathbb{R}^N)\hookrightarrow
L^q(\mathbb{R}^N)$$ is continuous. However $(V_2)$ together with $(V_1)$ implies that $E\hookrightarrow L^q(\mathbb{R}^N)$
is compact for all $q\geq 1$
(see Lemma 2.4 below). Surprisingly the assumption $(V_2)$ is much better than \\

\noindent $(V_2^\prime)$ $V(x)\ra+\infty$ as $|x|\ra+\infty$,\\

\noindent since $(V_1)$ together with $(V_2^\prime)$ only leads to
the compact embedding $E\hookrightarrow L^q(\mathbb{R}^N)$ for all
$q\geq N$ (see for example Costa \cite{Costa} for details). This is
the case in \cite{Adi-Yang,doo3}. However in this paper our argument
of proving main results seriously depends on the compact embedding
$E\hookrightarrow L^q(\mathbb{R}^N)$ for all $q\geq 1$.
\\

For any $\beta: 0\leq \beta<N$, we define a singular eigenvalue for
the $N$-Laplace operator by \be\label{lamda}\lambda_\beta=\inf_{u\in
E,\,u\not\equiv0}\f{\|u\|_E^N}
{\int_{\mathbb{R}^N}\f{|u|^N}{|x|^\beta}dx}.\ee It is easy to see
that
 $\lambda_\beta>0$. Write $m(r)=\sup_{|x|\leq r}V(x)$ and
  \be\label{m}\mathcal{M}=\inf_{r>0}\f{(N-\beta)^N}{\alpha_0^{N-1}r^{N-\beta}}e^{(N-\beta)m(r)\f{(N-2)!}{N^N}r^N},\ee
  where $\alpha_0$ is given by $(H_1)$. If $V(x)$ is continuous and $(V_1)$ is satisfied,
  then $m(r)$ is a positive continuous function and $\mathcal{M}$ can be attained by some $r>0$.\\

  One of our main results can be stated as follows:\\

 \noindent{\bf Theorem 1.1.} {\it Assume that $V(x)$ is a continuous function satisfying $(V_1)$ and $(V_2)$.
 $f:\mathbb{R}^N\times\mathbb{R}\ra\mathbb{R}$ is a continuous
 function and the hypothesis $(H_1)$, $(H_2)$ and $(H_3)$ hold. Furthermore we assume
 $$\leqno(H_4)\quad\quad\quad \limsup_{s\ra
0+}\f{N|F(x,s)|}{s^N}<\lambda_\beta\,\,\,{\rm
uniformly\,\,with\,\,respect\,\,to}\,\,\,x\in \mathbb{R}^N;$$
 $$\leqno(H_5)\quad\quad\quad \liminf_{s\ra
+\infty}sf(x,s)e^{-\alpha_0s^\f{N}{N-1}}=\beta_0>\mathcal{M}\,\,\,{\rm
uniformly\,\,with\,\,respect\,\,to}\,\,\,x\in \mathbb{R}^N.$$
 Then the equation (\ref{prob-0})
 has a nontrivial positive mountain-pass type weak solution.}\\

 Here and throughout this paper, we say that a weak solution $u$ is positive if $u(x)\geq 0$ for almost every $x\in\mathbb{R}^N$. It should be
 pointed out that $\mathcal{M}$ is not the best constant in $(H_5)$. It would be interesting if one can find
 an explicit smaller number replacing
 $\mathcal{M}$.\\

 In \cite{Adi-Yang}, Theorem A has been employed  to study a perturbation of the equation (\ref{prob-0}), namely

  \be\label{prob-1}
 -\Delta_Nu+V(x)|u|^{N-2}u=\f{f(x,u)}{|x|^\beta}+\epsilon h,\,\,\,x\in
  \mathbb{R}^N\,\,(N\geq 2),
\ee where $\epsilon>0$ is a constant and $h:\mathbb{R}^N\ra
\mathbb{R}$ is a function belonging to $E^*$, the dual space of $E$.
If $V(x)$ satisfies $(V_1)$, $(V_2^\prime)$, and $f(x,s)$ satisfies
$(H_1)-(H_4)$, then it was shown in \cite{Adi-Yang} that when
$\epsilon>0$ is sufficiently small and $h\not\equiv 0$, the problem
(\ref{prob-1}) has two weak solutions: one is of mountain-pass type
and the other is of negative energy. But we can not conclude that
the two solutions are distinct. In this paper, replacing
$(V_2^\prime)$ by $(V_2)$ and imposing additional condition $(H_5)$,
we can prove that the above two solutions are
distinct, namely \\

\noindent {\bf Theorem 1.2.} {\it  Suppose that $f(x,s)$ is
continuous in $\mathbb{R}^N\times\mathbb{R}$ and $(H_1)-(H_5)$ hold.
$V(x)$ is continuous in $\mathbb{R}^N$ satisfying $(V_1)$ and
$(V_2)$, $h$ belongs to $E^*$, the dual space of $E$, with $h\geq 0$
and $h\not\equiv 0$. Then there exists $\epsilon_0>0$ such that if
$0<\epsilon<\epsilon_0$, then the
problem (\ref{prob-1}) has two distinct positive weak solutions.}\\

The proof of Theorem 1.1 and Theorem 1.2 is based on Theorem A, the
mountain-pass theorem without the Palais-Smale condition
\cite{Rabin} and the Ekeland's variational principle \cite{Mawin},
which were also used in \cite{Adi-Yang,doo3}. Let us make some
reduction on problems (\ref{prob-0}) and (\ref{prob-1}). Set
$$\widetilde{f}(x,s)=\le\{\begin{array}{lll}0,\quad &f(x,s)<0\\[1.5ex]
f(x,s),& f(x,s)\geq 0.\end{array}\ri.$$ Assume $u\in E$ is a weak
solution of
\be\label{tild}-\Delta_Nu+V(x)|u|^{N-2}u=\f{\widetilde{f}(x,u)}{|x|^\beta}+\epsilon
h,\ee where $h\geq 0$ and $\epsilon> 0$, then the negative part of
$u$, namely
$$u_-(x)=\le\{\begin{array}{lll}0,\quad &u(x)>0\\[1.5ex]
u(x),& u(x)\leq 0\end{array}\ri.$$ belongs to the function space $E$
and satisfies \bna \int_{\mathbb{R}^N}(|\nabla
u_-|^N+V(x)|u_-|^N)dx&=&\int_{\mathbb{R}^N}\f{\widetilde{f}(x,u)}{|x|^\beta}u_-dx
+\epsilon\int_{\mathbb{R}^N}hu_-dx\\&=&\epsilon\int_{\mathbb{R}^N}hu_-dx\leq
0. \ena Hence $u_-(x)= 0$ for almost every $x\in \mathbb{R}^N$ and
thus $u$ is a positive weak solution of (\ref{tild}). This together with $(H_2)$ implies
$f(x,u)\geq 0$. It follows that $\widetilde{f}(x,u)=f(x,u)$.
Therefore $u$ is also a positive weak solution of (\ref{prob-1}).
When $h=0$, (\ref{prob-1}) becomes (\ref{prob-0}). Based on this, to
prove Theorems 1.1 and 1.2, it suffices to find weak solutions of
(\ref{prob-0}) and (\ref{prob-1}) with $f$ replaced by $\widetilde{f}$ respectively. So throughout this
paper, we can assume without loss of generality
\be\label{fzer0}f(x,s)\equiv 0,\quad\forall s<0.\ee

Before ending this introduction, we would like to mention that
results similar to Theorem 1.2 in two dimensional case, i.e. $N=2$,
was obtained by J. M. do \'O \cite{doo3}. Similar problems for
bi-Laplace equation in $\mathbb{R}^4$ was considered by the author
in \cite{Yang}. For compact Riemannian manifold case, we refer the
reader to \cite{doYang,Yangzhao}. Also it should be remarked that
results obtained in \cite{Adi-Yang} and in the present paper still
hold if there is only the subcritical case of (\ref{criticalineq}),
namely for any $\alpha<(1-\beta/N)\alpha_N$ and $\tau>0$,
$$\sup_{\int_{\mathbb{R}^N}(|\nabla
  u|^N+\tau|u|^N)dx\leq 1}
 \int_{\mathbb{R}^N}\f{e^{\alpha|u|^{N/(N-1)}}-\sum_{k=0}^{N-2}
 \f{\alpha^k|u|\,^{kN/(N-1)}}{k!}}{|x|^\beta}dx<\infty.$$
In fact, in \cite{Alves,Cao,doo1,Panda}, all the contributors only
used the above subcritical inequality.
 \\

The remaining part of this paper is organized as follows: In section
2, we display several key estimates in later compactness analysis.
In section 3, we consider the functionals related to problems
(\ref{prob-0}) and (\ref{prob-1}). Finally Theorem 1.1 is proved in
section 4 and Theorem 1.2 is proved in section 5.

\section{Key estimates}
In this section we will derive several technical lemmas for our use
later. For any integer $N\geq 2$ and real number $s$, we define a
function $\zeta:\mathbb{N}\times\mathbb{R}\ra \mathbb{R}$ by
\be\label{MF}\zeta(N,s)=e^s-\sum_{k=0}^{N-2}\f{s^k}{k!}=\sum_{k=N-1}^{\infty}\f{s^k}{k!}.\ee

\noindent{\bf Lemma 2.1.} {\it Let $s\geq 0$, $p\geq 1$ be real
numbers and $N\geq 2$ be an integer. Then there holds
\be\label{in}\le(\zeta(N,s)\ri)^p\leq \zeta(N,ps).\ee} {\it Proof.}
We prove (\ref{in}) by induction with respect to $N$. Define a
function
$$\phi(s)=(e^s-1)^p-(e^{ps}-1).$$
It is easy to see that for $s\geq 0$ and $p\geq 1$,
$$\phi^\prime(s)=p(e^s-1)^{p-1}-pe^{ps}\leq 0.$$
Hence $\phi(s)\leq\phi(0)=0$ and thus (\ref{in}) holds for $N=2$.
Suppose (\ref{in}) holds for $N\geq 2$, we only need to prove that
\be\label{intr}\le(\zeta(N+1,s)\ri)^p\leq \zeta(N+1,ps).\ee For this
purpose we set
$$\psi(s)=\le(e^s-\sum_{k=0}^{N-1}\f{s^k}{k!}\ri)^p-\le(e^{ps}-\sum_{k=0}^{N-1}\f{{(ps)}^k}{k!}\ri).$$
A straightforward calculation shows
\bna
 \psi^\prime(s)&=&p\le(e^s-\sum_{k=0}^{N-1}\f{s^k}{k!}\ri)^{p-1}\le(e^s-\sum_{k=1}^{N-1}\f{s^{k-1}}{(k-1)!}\ri)\\
 &&\quad-\le(p e^{ps}-p\sum_{k=1}^{N-1}\f{(ps)^{k-1}}{(k-1)!}\ri)\\
 &\leq& p\le\{\le(e^s-\sum_{k=1}^{N-1}\f{s^{k-1}}{(k-1)!}\ri)^p-\le(e^{ps}-\sum_{k=1}^{N-1}\f{(ps)^{k-1}}{(k-1)!}\ri)\ri\}\\
 &=&p\le\{\le(e^s-\sum_{k=0}^{N-2}\f{s^{k}}{k!}\ri)^p-\le(e^{ps}-\sum_{k=0}^{N-2}\f{(ps)^{k}}{k!}\ri)\ri\}\leq 0.
\ena
Here we have used the induction assumption $(\zeta(N,s))^p\leq \zeta(N,ps)$.
Thus $\psi(s)\leq \psi(0)=0$ for $s\geq 0$, and whence (\ref{intr}) holds. Therefore
(\ref{in}) holds for any integer $N\geq 2$. $\hfill\Box$ \\

\noindent{\bf Lemma 2.2.} {\it For all $N\geq 2$, $s\geq 0$, $t\geq
0$, $\mu>1$ and $\nu>1$ with $1/\mu+1/\nu=1$, there holds
$$\zeta(N,s+t)\leq \f{1}{\mu}\zeta(N,\mu s)+\f{1}{\nu}\zeta(N,\nu t).$$}
{\it Proof}. Observing that
$$\f{\p^2}{\p s^2}\zeta(2,s)=e^s\geq 0,\quad \f{\p^2}{\p s^2}\zeta(3,s)=e^s\geq 0$$
and when $N\geq 4$,
$$\f{\p^2}{\p s^2}\zeta(N,s)=e^s-\sum_{k=2}^{N-2}\f{s^{k-2}}{(k-2)!}=e^s-\sum_{k=0}^{N-4}\f{s^{k}}{k!}\geq 0,$$
we conclude that $\zeta(N,s)$ is convex with respect to $s$ for all $N\geq 2$. Hence
$$\zeta(N,s+t)=\zeta\le(N,\f{1}{\mu}\mu s+\f{1}{\nu}\nu t\ri)\leq \f{1}{\mu}\zeta(N,\mu s)+
\f{1}{\nu}\zeta(N,\nu t).$$
This concludes the lemma. $\hfill\Box$\\

\noindent{\bf Lemma 2.3.} {\it Let $(w_n)$ be a sequence in $E$.
Suppose $\|w_n\|_E=1$,
 $w_n\rightharpoonup w_0$ weakly in $E$, $w_n(x)\ra w_0(x)$ and $\nabla w_n(x)\ra \nabla w_0(x)$ for almost every
 $x\in\mathbb{R}^N$. Then for any $p: 0<p<\f{1}{(1-\|w_0\|_E^N)^{{1}/{(N-1)}}}$
 \be\label{2.1}\sup_n\int_{\mathbb{R}^N}\f{\zeta\le(N,\alpha_N(1-\beta/N)p|w_n|^{\f{N}{N-1}}\ri)}{|x|^\beta}dx<\infty.\ee}

\noindent{\it Proof.} Noticing that
$$|w_n|^{\f{N}{N-1}}=|w_n-w_0+w_0|^{\f{N}{N-1}}\leq (1+\epsilon)|w_n-w_0|^{\f{N}{N-1}}+c(\epsilon)|w_0|^{\f{N}{N-1}},$$
we have by using Lemma 2.2 \bna
\zeta\le(N,\alpha_N(1-\beta/N)p|w_n|^{\f{N}{N-1}}\ri)&\leq&
\f{1}{\mu}\zeta\le(N,\mu(1+\epsilon)\alpha_N(1-\beta/N)p|w_n-w_0|^{\f{N}{N-1}}\ri)
\\&&+\f{1}{\nu}\zeta\le(N,\nu c(\epsilon)\alpha_N(1-\beta/N)p|w_0|^{\f{N}{N-1}}\ri)\\
&\leq&\zeta\le(N,\mu(1+\epsilon)\alpha_N(1-\beta/N)p\|w_n-w_0\|_E^{\f{N}{N-1}}\le(\f{|w_n-w_0|}{\|w_n-w_0\|_E}\ri)^{\f{N}{N-1}}\ri)\\
&&+\zeta\le(N,\nu
c(\epsilon)\alpha_N(1-\beta/N)p|w_0|^{\f{N}{N-1}}\ri), \ena where
$\mu>1$, $\nu>1$ and $1/\mu+1/\nu=1$. By Br\'ezis-Lieb's Lemma
\cite{BL},
$$\|w_n-w_0\|_E^N=1-\|w_0\|_E^N+o_n(1),$$
where $o_n(1)\ra 0$ as $n\ra\infty$.
Hence for any $p: 0<p<\f{1}{(1-\|w_0\|_E^N)^{{1}/{(N-1)}}}$, one can
choose  $\epsilon>0$ sufficiently small and $\mu>1$ sufficiently
close to $1$ such that
$$\mu(1+\epsilon)\alpha_N(1-\beta/N)p\|w_n-w_0\|_E^{\f{N}{N-1}}<\alpha_N(1-\beta/N).$$
Now (\ref{2.1}) follows from Theorem A immediately. $\hfill\Box$\\

\noindent{\bf Lemma 2.4.} {\it Assume
$V:\mathbb{R}^N\times\mathbb{R}\ra \mathbb{R}$ is continuous and
$(V_1)$, $(V_2)$ hold.
 Then
$E$ is compactly embedded in $L^q(\mathbb{R}^N)$ for all $q\geq
1$.}\\

\noindent{\it Proof.} By $(V_1)$, the standard Sobolev embedding
  theorem implies that the following  embedding is continuous
  $$E\hookrightarrow W^{1,N}(\mathbb{R}^N)\hookrightarrow L^q(\mathbb{R}^N)\quad
  {\rm for\,\,all}\quad N\leq q<\infty.$$
  It follows from the H\"older inequality and $(V_2)$ that
  $$\int_{\mathbb{R}^N}|u|dx\leq\le(\int_{\mathbb{R}^N}\f{1}{V^{\f{1}{N-1}}}dx\ri)^{1-1/N}
  \le(\int_{\mathbb{R}^N}V|u|^Ndx\ri)^{1/N}\leq \le(\int_{\mathbb{R}^N}\f{1}{V^{\f{1}{N-1}}}dx\ri)^{1-1/N}
  \|u\|_E.$$
  For any $\gamma:1<\gamma<N$, there holds
  $$\int_{\mathbb{R}^N}|u|^\gamma dx\leq \int_{\mathbb{R}^N}(|u|+|u|^N)dx\leq \le(\int_{\mathbb{R}^N}\f{1}{V^{\f{1}{N-1}}}
  dx\ri)^{1-1/N}\|u\|_E
  +\f{1}{V_0}\|u\|_E^N,$$
  where $V_0$ is given by $(V_1)$.
  Thus we get continuous embedding $E\hookrightarrow
  L^q(\mathbb{R}^N)$ for all $q\geq 1$.

  To prove that the above embedding is also compact, take a sequence of functions
  $(u_k)\subset E$ such that $\|u_k\|_E\leq
  C$ for all $k$, we must prove that up to a subsequence there exists some $u\in E$ such that $u_k$
  convergent to  $u$ strongly in $L^q(\mathbb{R}^N)$
  for all $q\geq 1$. Without loss of generality we may assume
  \be\label{bdd}\le\{\begin{array}{lll}
  u_k\rightharpoonup u\quad&{\rm weakly\,\,in}\quad E\\[1.5ex]
  u_k\rightarrow u\quad&{\rm strongly\,\,in}\quad L^q_{\rm loc}(\mathbb{R}^N),\,\,\forall q\geq 1\\[1.5ex]
  u_k\ra u\quad&{\rm almost\,\,everywhere\,\,in}\quad \mathbb{R}^N.
  \end{array}\ri.
  \ee
   In view of $(V_2)$, for any $\epsilon>0$, there exists $R>0$ such that
  $$\le(\int_{|x|> R}\f{1}{V^{\f{1}{N-1}}}dx\ri)^{1-1/N}<\epsilon.$$
  Hence
  \be\label{res}\int_{|x|>R}|u_k-u|dx\leq \le(\int_{\mathbb{R}^N}\f{1}{V^{\f{1}{N-1}}}dx\ri)^{1-1/N}
  \le(\int_{\mathbb{R}^N}V|u|^Ndx\ri)^{1/N}\leq \epsilon\|u_k-u\|_E\leq
  C\epsilon.\ee
  Here and in the sequel we often denote various constants by the same $C$.
  On the other hand, it follows from (\ref{bdd}) that $u_k\ra u$ strongly in $L^1(\mathbb{B}_R(0))$,
  where $\mathbb{B}_R(0)\subset\mathbb{R}^N$ is the ball centered at $0$ with radius $R$.
  This together with (\ref{res}) leads to
  $$\limsup_{k\ra \infty}\int_{\mathbb{R}^N}|u_k-u|dx\leq C\epsilon.$$
  Since $\epsilon$ is arbitrary, we obtain
  $$\lim_{k\ra \infty}\int_{\mathbb{R}^N}|u_k-u|dx=0.$$
  For $q>1$, it follows from the continuous embedding $E\hookrightarrow L^s(\mathbb{R}^N)$ ($s\geq 1$)
  that
  \bna
  \int_{\mathbb{R}^N}|u_k-u|^qdx&=&\int_{\mathbb{R}^N}|u_k-u|^\f{1}{2}|u_k-u|^{(q-\f{1}{2})}dx\\
  &\leq&\le(\int_{\mathbb{R}^N}|u_k-u|dx\ri)^{1/2}\le(\int_{\mathbb{R}^N}|u_k-u|^{2q-1}dx\ri)^{1/2}\\
  &\leq&C\le(\int_{\mathbb{R}^N}|u_k-u|dx\ri)^{1/2}\ra 0
  \ena
  as $k\ra\infty$. This concludes the lemma. $\hfill\Box$
  \\

\section{Functionals and compactness analysis}

\subsection{The functionals and their profiles\\}

  As we mentioned in the introduction, problems (\ref{prob-0}) and (\ref{prob-1})
  have variational structure. To apply the critical point theory, we define the functional $J_{\beta,\,\epsilon}: E\ra \mathbb{R}$ by
  $$J_{\beta,\,\epsilon}(u)=\f{1}{N}\|u\|_E^N-\int_{\mathbb{R}^N}\f{F(x,u)}{|x|^\beta}dx-\epsilon\int_{\mathbb{R}^N}hudx,$$
  where $\epsilon\geq 0$, $0\leq\beta<N$, $\|u\|_E$ is the norm of $u\in E$ defined by (\ref{E-norm}) and $F(x,s)=\int_0^sf(x,t)dt$ is the
  primitive of $f(x,s)$.
  Assume $f$ satisfies the hypothesis $(H_1)$. Then there exist some positive constants $\alpha_1>\alpha_0$ and
  $b_3$ such that for all $(x,s)\in \mathbb{R}^N\times\mathbb{R}$,
  $F(x,s)\leq
  b_3\zeta(N,\alpha_1|s|^{N/(N-1)})$.
  Thus $J_{\beta,\,\epsilon}$ is well defined thanks to Theorem A. In the case $\epsilon=0$, we denote
  $J_{\beta,0}$ for simplicity by
   $$J(u)=\f{1}{N}\|u\|_E^N-\int_{\mathbb{R}^N}\f{F(x,u)}{|x|^\beta}dx.$$
   The profiles of the functionals $J_{\beta,\epsilon}$ and $J(u)$
  are well described in the following lemma.  \\

  \noindent{\bf Lemma 3.1.} {\it Assume $(V_1)$, $(H_1)$, $(H_2)$, $(H_3)$ and
  $(H_4)$ are satisfied. Then $(i)$ for any nonnegative, compactly supported function $u\in
  W^{1,N}(\mathbb{R}^N)\setminus\{0\}$, there holds $J_{\beta,\,\epsilon}(tu)\ra
  -\infty$ as $t\ra +\infty$; $(ii)$ there exists $\epsilon_1>0$
 such that when $0<\epsilon<\epsilon_1$, one can find
  $r_\epsilon$, $\vartheta_\epsilon>0$ such that
   $J_{\beta,\,\epsilon}(u)\geq\vartheta_\epsilon$ for all $u$ with $\|u\|_E=r_\epsilon$,
   where $r_\epsilon$ can be further chosen such that $r_\epsilon\ra 0$ as $\epsilon\ra
   0$. When $\epsilon=0$, there exists $\delta>0$ and $r>0$ such that $J(u)\geq \delta$ for
   all $\|u\|_E=r$; $(iii)$ assume $\epsilon>0$ and $h\not\equiv 0$, there exists a constant $\tau>0$ such that if
   $0<t<\tau$, then
   $\inf_{\|u\|_{E}\leq t}J_{\beta,\epsilon}(u)<0$.}\\

   \noindent{\it Proof}.  We refer the reader to (\cite{Adi-Yang}, Lemmas 4.1, 4.2 and 4.3) for details.
   It is remarkable that we can also apply
   Lemma 2.1 and Lemma 2.4 instead of decreasing rearrangement argument
   in the proof of (\cite{Adi-Yang}, Lemma 4.2) and thus simplify it. $\hfill\Box$\\

  \noindent To use the critical point theory, we need some regularity of the functionals $J_{\beta,\epsilon}$
  and $J$. In fact, by Proposition 1 in \cite{doo3} and
  standard arguments (see for example \cite{Rabin}), one can see that both $J_{\beta,\epsilon}$ and $J$
  belong to $\mathcal{C}^1(E,\mathbb{R})$.
  A straightforward calculation shows
   \bea\label{j0}
    &&\langle J^\prime(u),\phi\rangle=\int_{\mathbb{R}^N}\le(|\nabla
    u|^{N-2}\nabla u\nabla
    \phi+V|u|^{N-2}u\phi\ri)dx-\int_{\mathbb{R}^N}\f{f(x,u)}{|x|^{\beta}}\phi d
    x,\\\label{j1}&&\langle J_{\beta,\,\epsilon}^\prime(u),\phi\rangle=\int_{\mathbb{R}^N}\le(|\nabla
    u|^{N-2}\nabla u\nabla
    \phi+V|u|^{N-2}u\phi\ri)dx-\int_{\mathbb{R}^N}\f{f(x,u)}{|x|^{\beta}}\phi d
    x-\epsilon\int_{\mathbb{R}^N} h\phi dx\qquad
   \eea
   for all $\phi\in E$. Hence weak solutions of (\ref{prob-0}) and (\ref{prob-1}) are
   critical points of $J$ and $J_{\beta,\,\epsilon}$ respectively.\\

  \subsection{Min-Max value\\}

  In this subsection, we prepare for estimating the min-max value of the functionals $J$ and $J_{\beta,\epsilon}$.
  The idea is to construct a sequence of functions $M_n\in E$ and estimate $\max\limits_{t\geq 0} J(t M_n)$ and
  $\max\limits_{t\geq 0} J_{\beta,\epsilon}(tM_n)$.
  Recall Moser's function sequence
  $$\widetilde{M}_n(x,r)=\f{1}{\omega_{N-1}^{1/N}}\le\{
  \begin{array}{lll}
  (\log n)^{1-1/N},\quad &|x|\leq r/n\\[1.5ex]
  \f{\log \f{r}{|x\,|}}{(\log n)^{1/N}}, &r/n<|x|\leq r\\[1.5ex]
  0,&|x|>r.
  \end{array}
  \ri.$$
  Let $M_n(x,r)=\f{1}{\|\widetilde{M}_n\|_E}\widetilde{M}_n(x,r)$.
  Then $M_n$ belongs to $E$ with its support in $\mathbb{B}_{r}(0)$ and $\|M_n\|_E=1$.\\

  \noindent{\bf Lemma 3.2.} {\it Assume $V(x)$ is continuous and $(V_1)$ is satisfied. Then there holds
  \be
  \label{nor}\|\widetilde{M}_n\|_E^N\leq 1+\f{m(r)}{\log n}\le(\f{(N-1)!}{N^N}r^N+o_n(1)\ri),
  \ee
  where $m(r)=\max\limits_{|x|\leq r}V(x)$ and $o_n(1)\ra 0$ as $n\ra\infty$.}\\

  \noindent{\it Proof.} It is easy to calculate
  $$\int_{\mathbb{R}^N}|\nabla\widetilde{M}_n|^Ndx=\f{1}{\omega_{N-1}}\int_{\f{r}{n}\leq |x|\leq r}
  \f{1}{|x|^N\log n}dx=1.$$
  Integration by parts gives
  \bna
  \int_{\f{r}{n}\leq |x|\leq r}\le(\log\f{r}{|x|}\ri)^Ndx&=&\omega_{N-1}\int_{\f{r}{n}}^rs^{N-1}\le(\log\f{r}{s}\ri)^Nds\\
  &=&-\f{\omega_{N-1}}{N}\le(\f{r}{n}\ri)^N\le(\log n\ri)^N+\omega_{N-1}\int_{\f{r}{n}}^rs^{N-1}\le(\log\f{r}{s}\ri)^{N-1}ds\\
  &=&-{\omega_{N-1}}\le(\f{r}{n}\ri)^N\le\{\f{1}{N}(\log n)^N+\f{1}{N}(\log n)^{N-1}+\f{N-1}{N^2}(\log n)^{N-2}\ri.\\
  &&\quad\quad\quad\quad\quad\le.+\cdots+\f{(N-1)(N-2)\cdots 3}{N^{N-2}}(\log n)^2\ri\}\\
  &&+\omega_{N-1}\f{(N-1)!}{N^{N-2}}\int_{\f{r}{n}}^rs^{N-1}\log\f{r}{s}ds\\
  &=&\omega_{N-1}\f{(N-1)!}{N^N}r^N+o_n(1).
  \ena
  Hence
  \bna
  \int_{\mathbb{R}^N}|\widetilde{M}_n|^Ndx&=&\f{1}{\omega_{N-1}}\int_{|x|\leq r/n}(\log n)^{N-1}dx+
  \f{1}{\omega_{N-1}}\int_{\f{r}{n}\leq |x|\leq r}\f{\le(\log\f{r}{|x|}\ri)^N}{\log n}dx\\
  &=&\le(\f{r}{n}\ri)^N\f{(\log n)^{N-1}}{N}+\f{1}{\omega_{N-1}\log n}\int_{\f{r}{n}\leq |x|\leq r}\le(\log\f{r}{|x|}\ri)^Ndx\\
  &=&\f{1}{\log n}\le(\f{(N-1)!}{N^N}r^N+o_n(1)\ri),
  \ena
  and thus
  \bna
  \|\widetilde{M}_n\|_E^N&=&\int_{\mathbb{R}^N}|\nabla \widetilde{M}_n|^Ndx+\int_{\mathbb{R}^N}V(x)|\widetilde{M}_n|^Ndx\\
  &\leq& 1+m(r)\int_{\mathbb{R}^N}|\widetilde{M}_n|^Ndx\\
  &=&1+\f{m(r)}{\log n}\le(\f{(N-1)!}{N^N}r^N+o_n(1)\ri).
  \ena
  This is exactly (\ref{nor}). $\hfill\Box$\\

  \noindent{\bf Lemma 3.3.} {\it Assume $(V_1)$, $(H_1)$, $(H_2)$, $(H_3)$ and $(H_5)$. There exists some $n\in\mathbb{N}$ such that
  \be\label{mx}
  \max_{t\geq 0}J(tM_n)<\f{1}{N}\le(\f{N-\beta}{N}\f{\alpha_N}{\alpha_0}\ri)^{N-1}.
  \ee
  Furthermore for the above $n$ there exists some $\epsilon^*>0$ and $\delta^*>0$
  such that if $0\leq \epsilon<\epsilon^*$, then
  \be\label{mx-eps}
  \max_{t\geq 0}J_{\beta,\epsilon}(tM_n)<\f{1}{N}\le(\f{N-\beta}{N}\f{\alpha_N}{\alpha_0}\ri)^{N-1}-\delta^*.
  \ee
  }
  {\it Proof.} We first prove (\ref{mx}).
  By $(H_5)$ and (\ref{m}) (the definition of $\mathcal{M}$), there exists some $r>0$ such that
  \be\label{beta0}\beta_0>\f{(N-\beta)^N}{\alpha_0^{N-1}r^{N-\beta}}e^{(N-\beta)m(r)\f{(N-2)!}{N^N}r^N}.\ee
  Suppose by contradiction that for all $n\in\mathbb{N}$
  \be\label{cn}\max_{t\geq 0}J(tM_n)\geq\f{1}{N}\le(\f{N-\beta}{N}\f{\alpha_N}{\alpha_0}\ri)^{N-1}.\ee
  By $(i)$ of Lemma 3.1, $\forall n\in\mathbb{N}$, there exists $t_n>0$ such that
  $$J(t_nM_n)=\max_{t\geq 0}J(tM_n).$$
  Thus (\ref{cn}) gives
  $$J(t_nM_n)=\f{t_n^N}{N}-\int_{\mathbb{R}^N}\f{F(x,t_nM_n)}{|x|^\beta}dx\geq \f{1}{N}
  \le(\f{N-\beta}{N}\f{\alpha_N}{\alpha_0}\ri)^{N-1}.$$
  Noticing that $F(x,\cdot)\geq 0$, we have
  \be\label{gg}t_n^N\geq \le(\f{N-\beta}{N}\f{\alpha_N}{\alpha_0}\ri)^{N-1}.\ee
  It is easy to see that at $t=t_n$,
  $$\f{d}{dt}\le(\f{t^N}{N}-\int_{\mathbb{R}^N}\f{F(x,tM_n)}{|x|^\beta}dx\ri)=0,$$
  or equivalently
  \be\label{der}t_n^N=\int_{\mathbb{R}^N}\f{t_nM_nf(x,t_nM_n)}{|x|^\beta}dx.\ee
  By $(H_5)$, $\forall \eta>0$, $\exists R_\eta>0$ such that for all $x\in \mathbb{R}^N$ and $u\geq R_\eta$
  \be\label{rep}uf(x,u)\geq (\beta_0-\eta)e^{\alpha_0 |u|^{\f{N}{N-1}}}.\ee
  By Lemma 3.2, when $|x|\leq \f{r}{n}$, we have
  \bea\nonumber
  M_n^{\f{N}{N-1}}(x,r)&\geq&\f{1}{\omega_{N-1}^{\f{1}{N-1}}}\f{\log n}{1+\f{1}{N-1}\f{m(r)}{\log n}\le(\f{(N-1)!}{N^N}r^N+o_n(1)\ri)}\\
  \label{nr}&=&\omega_{N-1}^{-\f{1}{N-1}}\log n-\omega_{N-1}^{-\f{1}{N-1}}{m(r)}\f{(N-2)!}{N^N}r^N+o_n(1).
  \eea
  Hence we have by combining (\ref{der}) and (\ref{rep}) that
  \bea{\nonumber}
  t_n^N&\geq&(\beta_0-\eta)\int_{|x|\leq \f{r}{n}}\f{e^{\alpha_0|t_nM_n|^{\f{N}{N-1}}}}{|x|^\beta}dx\\{\nonumber}
  &=&(\beta_0-\eta)\int_{|x|\leq \f{r}{n}}\f{e^{\alpha_0\omega_{N-1}^{-\f{1}{N-1}}t_n^{\f{N}{N-1}}
  \le(\log n-m(r)\f{(N-2)!}{N^N}r^N+o_n(1)\ri)}}{|x|^\beta}dx\\\label{3}
  &=&(\beta_0-\eta)\f{\omega_{N-1}}{N-\beta}\le(\f{r}{n}\ri)^{N-\beta}
  e^{\alpha_0\omega_{N-1}^{-\f{1}{N-1}}t_n^{\f{N}{N-1}}
  \le(\log n-m(r)\f{(N-2)!}{N^N}r^N+o_n(1)\ri)}.
  \eea
  This yields that $t_n$ is a bounded sequence. In view of (\ref{gg}), we can also see from (\ref{3}) that
  \be\label{lim}\lim_{n\ra\infty}t_n^N= \le(\f{N-\beta}{N}\f{\alpha_N}{\alpha_0}\ri)^{N-1}.\ee
  For otherwise there exists some $\delta>0$ such that for sufficiently large $n$
  $$t_n^N\geq \le(\delta+\f{N-\beta}{N}\f{\alpha_N}{\alpha_0}\ri)^{N-1}.$$
  Thus
  $$\alpha_0\omega_{N-1}^{-\f{1}{N-1}}t_n^{\f{N}{N-1}}\geq N-\beta+\alpha_0\omega_{N-1}^{-\f{1}{N-1}}\delta$$
  and whence the right hand of (\ref{3}) tends to infinity which contradicts the bounded-ness of $t_n$.

  Now we estimate $\beta_0$. It follows from (\ref{der}) and (\ref{rep}) that
  \bea
  t_n^{N}&\geq& (\beta_0-\eta)\int_{|x|\leq r}\f{e^{\alpha_0|t_nM_n|^{\f{N}{N-1}}}}{|x|^\beta}dx+
  \int_{t_nM_n<R_\eta}\f{t_nM_nf(x,t_nM_n)}{|x|^\beta}dx{\nonumber}\\\label{tn}&&\qquad
  -(\beta_0-\eta)\int_{t_nM_n<R_\eta}\f{e^{\alpha_0|t_nM_n|^{\f{N}{N-1}}}}{|x|^\beta}dx.
  \eea
  Since $M_n\ra 0$ almost everywhere in $\mathbb{R}^N$, we have by using the Lebesgue's dominated convergence theorem
  \bea\label{p1}&&\lim_{n\ra\infty}\int_{t_nM_n<R_\eta}\f{t_nM_nf(x,t_nM_n)}{|x|^\beta}dx=0,\\&&
  \label{p2}\lim_{n\ra\infty}\int_{t_nM_n<R_\eta}\f{e^{\alpha_0|t_nM_n|^{\f{N}{N-1}}}}{|x|^\beta}dx=
  \int_{|x|\leq r}\f{1}{|x|^\beta}dx=\f{\omega_{N-1}r^{N-\beta}}{N-\beta}.\eea
  Using (\ref{gg}),
  \be\label{p3}
  \int_{|x|\leq r}\f{e^{\alpha_0|t_nM_n|^{\f{N}{N-1}}}}{|x|^\beta}dx\geq \int_{|x|\leq \f{r}{n}}
  \f{e^{\alpha_N(1-\beta/N)M_n^{\f{N}{N-1}}}}{|x|^\beta}dx+\int_{\f{r}{n}\leq|x|\leq {r}}
  \f{e^{\alpha_N(1-\beta/N)M_n^{\f{N}{N-1}}}}{|x|^\beta}dx.
  \ee
  On one hand we have by (\ref{nr})
  \bna
  \int_{|x|\leq \f{r}{n}}
  \f{e^{\alpha_N(1-\beta/N)M_n^{\f{N}{N-1}}}}{|x|^\beta}dx&\geq& e^{\alpha_N(1-\beta/N)\omega_{N-1}^{-\f{1}{N-1}}\log n-\omega_{N-1}^{-\f{1}{N-1}}{m(r)}\f{(N-2)!}{N^N}r^N+o_n(1)}\int_{|x|\leq \f{r}{n}}\f{1}{|x|^\beta}dx\\
  &=&\f{\omega_{N-1}}{N-\beta}\le(\f{r}{n}\ri)^{N-\beta}e^{(N-\beta)\log n-(N-\beta)m(r)\f{(N-2)!}{N^N}r^N+o_n(1)}\\
  &=&\f{\omega_{N-1}r^{N-\beta}}{N-\beta}e^{-(N-\beta)m(r)\f{(N-2)!}{N^N}r^N+o_n(1)}.
  \ena
  On the other hand, by definition of $M_n$,
  \bna
  \int_{\f{r}{n}\leq|x|\leq {r}}
  \f{e^{\alpha_N(1-\beta/N)M_n^{\f{N}{N-1}}}}{|x|^\beta}dx&=&
  \int_{\f{r}{n}\leq|x|\leq {r}}
  \f{e^{(N-\beta)\le((\log n)^{-1/N}\|\widetilde{M}_n\|_E^{-1}\log \f{r}{|x|}\ri)^{\f{N}{N-1}}}}{|x|^\beta}dx\\
  &=&\omega_{N-1}\int_{\f{r}{n}}^rt^{N-\beta-1}e^{(N-\beta)\le((\log n)^{-1/N}\|\widetilde{M}_n\|_E^{-1}\log \f{r}{t}\ri)^{\f{N}{N-1}}}dt\\
  &=&\omega_{N-1}r^{N-\beta}\int_0^{{(\log n)^{1-1/N}\|\widetilde{M}_n\|_E^{-1}}}(\log n)^{1/N}\|\widetilde{M}_n\|_E\\
  &&\qquad e^{(N-\beta)s^{\f{N}{N-1}}
  -(N-\beta)\|\widetilde{M}_n\|_E(\log n)^{1/N}s}ds\\
  &\geq&\omega_{N-1}r^{N-\beta}\int_0^{{(\log n)^{1-1/N}\|\widetilde{M}_n\|_E^{-1}}}(\log n)^{1/N}\|\widetilde{M}_n\|_E\\
  &&\qquad e^{
  -(N-\beta)\|\widetilde{M}_n\|_E(\log n)^{1/N}s}ds\\
  &=&\f{\omega_{N-1}r^{N-\beta}}{N-\beta}\le(1-e^{-(N-\beta)\log n}\ri).
  \ena
  Here we have used the change of variable $t=re^{-\|\widetilde{M}_n\|_E(\log n)^{1/N}s}$ in the third equality. Hence we obtain by passing to the limit
  $n\ra\infty$ in (\ref{p3})
  $$\liminf_{n\ra\infty}\int_{|x|\leq r}\f{e^{\alpha_0|t_nM_n|^{\f{N}{N-1}}}}{|x|^\beta}dx\geq
  \f{\omega_{N-1}r^{N-\beta}}{N-\beta}\le(1+e^{-(N-\beta)m(r)\f{(N-2)!}{N^N}r^N}\ri).$$
  This together with (\ref{lim})-(\ref{p2}) implies
  $$\le(\f{N-\beta}{N}\f{\alpha_N}{\alpha_0}\ri)^{N-1}\geq (\beta_0-\eta)\f{\omega_{N-1}r^{N-\beta}}{N-\beta}e^{-(N-\beta)m(r)\f{(N-2)!}{N^N}r^N}.$$
  Since $\eta>0$ is arbitrary, we have
  $$\beta_0\leq \f{(N-\beta)^N}{\alpha_0^{N-1}r^{N-\beta}}e^{(N-\beta)m(r)\f{(N-2)!}{N^N}r^N}.$$
  This contradicts (\ref{beta0}) and ends the proof of (\ref{mx}).

  Secondly it follows from (\ref{mx}) and the definition of $J_{\beta,\epsilon}$ that (\ref{mx-eps})
  holds. $\hfill\Box$\\

  \subsection{Palais-Smale sequence\\}

  In this subsection, we will show that the weak limit of a Palais-Smale sequence for $J_{\beta,\epsilon}$ is the weak solution
  of $(\ref{prob-1})$. (Respectively the weak limit of a Palais-Smale sequence for $J$ is also the weak solution
  of $(\ref{prob-0})$.)\\

  \noindent{\bf Lemma 3.4.} {\it Assume that $(V_1)$, $(V_2)$, $(H_1)$, $(H_2)$ and $(H_3)$ are satisfied.
  Let $(u_n)\subset E$ be an arbitrary Palais-Smale sequence of $J_{\beta,\epsilon}$, i.e.,
   \be\label{PS}J_{\beta,\epsilon}(u_n)\ra c,\,\, J^\prime_{\beta,\epsilon}(u_n)\ra 0 \,\, {\rm in}\,\,
   E^*
   \,\,{\rm as}\,\, n\ra \infty,\ee
   where $E^*$ denotes the dual space of $E$.
   Then there exist a subsequence of $(u_n)$ (still denoted by $(u_n)$) and $u\in E$ such that
   $u_n\rightharpoonup u$ weakly in $E$, $u_n\ra u$ strongly in
   $L^q(\mathbb{R}^N)$ for all $q\geq 1$, and
   \bna\le\{\begin{array}{lll}
   \nabla u_n(x)\ra\nabla u(x)\quad{\rm
   a.\,\,e.\,\,\,in}\quad\mathbb{R}^N\\[1.5ex]
   \f{f(x,\,u_n)}{|x|^\beta}\ra \f{f(x,\,u)}{|x|^\beta}\,\,{\rm strongly\,\,in}\,\,
   L^1(\mathbb{R}^N)\\[1.5ex]
    \f{F(x,\,u_n)}{|x|^\beta}\ra \f{F(x,\,u)}{|x|^\beta}\,\, {\rm strongly\,\,in}\,\,
   L^1(\mathbb{R}^N).
   \end{array}\ri.\ena
   Furthermore $u$ is a weak solution of (\ref{prob-1}). The same conclusion holds when $\epsilon=0$.}\\

   \noindent {\it Proof.} Assume $(u_n)$ is a Palais-Smale sequence of $J_{\beta,\epsilon}$. By $(\ref{PS})$, we have
   \bea\label{1}
   &&\f{1}{N}\|u_n\|_E^N-
   \int_{\mathbb{R}^N}\f{F(x,u_n)}{|x|^\beta}dx-\epsilon\int_{\mathbb{R}^N}hu_ndx\ra
   c\,\,{\rm as}\,\,n\ra\infty,\\\label{2}
   &&\le|\int_{\mathbb{R}^N}\le(|\nabla u_n|^{N-2}\nabla u_n\nabla\psi+V|u_n|^{N-2}u_n\psi\ri)dx
   -\int_{\mathbb{R}^N}\f{f(x,u_n)}{|x|^\beta}\psi
   dx-\epsilon\int_{\mathbb{R}^N}h\psi dx\ri|\leq \tau_n\|\psi\|_E\qquad\quad
   \eea
   for all $\psi\in E$, where $\tau_n\ra 0$ as $n\ra\infty$. Noticing that (\ref{fzer0}),
   we have by $(H_2)$ that $0\leq \mu F(x,u_n)\leq
   u_nf(x,u_n)$ for some $\mu>N$. Taking $\psi=u_n$ in (\ref{2}) and multiplying
   (\ref{1}) by $\mu$, we have
   \bna
   \le(\f{\mu}{N}-1\ri)\|u_n\|_E^N&\leq&\le(\f{\mu}{N}-1\ri)\|u_n\|_E^N-\int_{\mathbb{R}^N}\f{\mu F(x,u_n)-f(x,u_n)u_n}{|x|^\beta}dx\\
   &\leq&\mu |c|+\tau_n\|u_n\|_E+(\mu+1)\epsilon\|h\|_{E^*}\|u_n\|_E
   \ena
   Therefore $\|u_n\|_E$ is bounded. It then follows from  (\ref{1}), (\ref{2}) that
   \be\label{c}\int_{\mathbb{R}^N}\f{f(x,u_n)u_n}{|x|^\beta}
   dx\leq C,\quad \int_{\mathbb{R}^N}\f{F(x,u_n)}{|x|^\beta}dx\leq
   C\ee for some constant $C$ depending only on $\mu$, $N$ and $\|h\|_{E^\ast}$.
    By Lemma 2.4, up to a subsequence, $u_n\ra u$ strongly in $L^q(\mathbb{R}^N)$ for some $u\in E$, $\forall q\geq 1$.
    This immediately leads to $u_n\ra u$ almost everywhere in $\mathbb{R}^N$.
    Now we {\it claim} that up to a subsequence
   \be\label{f-conv}\lim_{n\ra\infty}\int_{\mathbb{R}^N}\f{|f(x,\,u_n)-f(x,\,u)|}{|x|^\beta}dx=0.\ee
   In fact, since $f(x,\cdot)\geq 0$, it suffices to prove that up to a subsequence
   \be\label{fc}\lim_{n\ra\infty}\int_{\mathbb{R}^N}\f{f(x,\,u_n)}{|x|^\beta}dx=
   \lim_{n\ra\infty}\int_{\mathbb{R}^N}\f{f(x,\,u)}{|x|^\beta}dx.\ee
   Since $u,\,\f{f(x,u)}{|x|^\beta}\in L^1(\mathbb{R}^N)$, we have
  $$\lim_{\eta\ra+\infty}\int_{|u|\geq\eta}\f{f(x,u)}{|x|^\beta}dx=0.$$
  Let $C$ be the constant in (\ref{c}). Given any $\delta>0$, one can select some $M>{C}/{\delta}$ such that
  \be\label{abs}
  \int_{|u|\geq M}\f{f(x,u)}{|x|^\beta}dx<\delta.
  \ee
  It follows from (\ref{c}) that
  \be\label{abe}\int_{|u_n|\geq M}\f{f(x,u_n)}{|x|^\beta}dx\leq \f{1}{M}\int_{|u_n|\geq
  M}\f{f(x,u_n)u_n}{|x|^\beta}dx<\delta.\ee
  For all $x\in\{x\in\mathbb{R}^N:|u_n|<M\}$, by our assumption $(H_1)$, there exists a constant $C_1$ depending only on
  $M$ such that $|f(x,u_n(x))|\leq
  C_1|u_n(x)|^{N-1}$. Notice that $|x|^{-\beta}{|u_n|^{N-1}}\ra |x|^{-\beta}|u|^{N-1}$ strongly in
  $L^{1}(\mathbb{R}^N)$ and $u_n\ra u$ almost everywhere
  in $\mathbb{R}^N$.
  By the generalized Lebesgue's dominated convergence theorem, we obtain
  \be\label{domin}\lim_{n\ra
  \infty}\int_{|u_n|<M}\f{f(x,u_n)}{|x|^\beta}dx=\int_{|u|<M}\f{f(x,u)}{|x|^\beta}dx.\ee
  Combining (\ref{abs}), (\ref{abe}) and (\ref{domin}), we can find
  some $K>0$ such that when $n>K$,
  $$\le|\int_{\mathbb{R}^N}\f{f(x,u_n)}{|x|^\beta}dx-\int_{\mathbb{R}^N}\f{f(x,u)}{|x|^\beta}dx\ri|<3\delta.$$
  Hence (\ref{fc}) holds and thus our claim (\ref{f-conv}) holds.
   By $(H_1)$ and $(H_3)$, there exist constants $c_1$, $c_2>0$
   such that
   $$F(x,u_n)\leq c_1|u_n|^N+c_2f(x,u_n).$$
   In view of (\ref{f-conv}) and Lemma 2.4, it follows from the generalized Lebesgue's dominated convergence
   theorem
   $$\lim_{n\ra\infty}\int_{\mathbb{R}^N}\f{|F(x,\,u_n)-F(x,\,u)|}{|x|^\beta}dx=0.$$
   Using the argument of proving (4.26) in \cite{Adi-Yang}, we have
   $\nabla u_n(x)\ra\nabla u(x)$ a. e. in $\mathbb{R}^N$
   and
   $$|\nabla u_n|^{N-2}\nabla u_n\rightharpoonup |\nabla u|^{N-2}\nabla u\quad{\rm weakly\,\,in}\quad \le(L^{\f{N}{N-1}}
   (\mathbb{R}^N)\ri)^N.$$
   Finally passing to the limit $n\ra \infty$ in $(\ref{2})$, we have
   $$\int_{\mathbb{R}^N}\le(|\nabla u|^{N-2}\nabla u\nabla\psi+V|u|^{N-2}u\psi\ri)dx-\int_{\mathbb{R}^N}\f{f(x,u)}{|x|^\beta}\psi
   dx-\epsilon\int_{\mathbb{R}^4}h\psi dx=0$$
   for all $\psi\in C_0^\infty(\mathbb{R}^N)$, which is dense in $E$. Hence $u$ is a weak solution of
   $(\ref{prob-1})$. After checking the above argument, $\epsilon$ need not to be nonzero, i.e. the
   same conclusion holds for $J$. $\hfill\Box$\\

   \noindent{\bf Remark 3.5.} Similar results of Lemma 3.4 was also
   established by J. M. do \'O in two
   dimensional case \cite{do-de} and by the author for bi-Laplace equation in four
   dimensional Euclidean space \cite{Yang}.

\section{Nontrivial positive solution\\}

In this section, we will prove Theorem 1.1.
 It suffices to look for nontrivial
 critical points of the functional $J$ in the function space $E$.\\

\noindent{\it Proof of Theorem 1.1.}
 By $(i)$ and $(ii)$ of Lemma 3.1, $J$ satisfies all the hypothesis of
 the mountain-pass theorem except for the Palais-Smale condition: $J\in
   \mathcal{C}^1(E,\mathbb{R})$; $J(0)=0$; $J(u)\geq \delta>0$ when
  $\|u\|_E=r$; $J(e)<0$ for some  $e\in E$  with
   $\|e\|_E>r$.
    Then using the mountain-pass theorem
   without the Palais-Smale condition \cite{Rabin}, we can find a sequence
   $(u_n)$ in $E$ such that
   $$J(u_n)\ra c>0,\quad J^\prime(u_n)\ra 0\,\,{\rm in}\,\, E^*,$$
   where
   $$c=\min_{\gamma\in\Gamma}\max_{u\in\gamma}J(u)\geq \delta$$
   is the min-max value of $J$, where $\Gamma=\{\gamma\in\mathcal{C}([0,1],E): \gamma(0)=0,
   \gamma(1)=e\}$. By (\ref{j0}), this is equivalent to saying
   \bea\label{01}
   &&\f{1}{N}\|u_n\|_E^N-
   \int_{\mathbb{R}^N}\f{F(x,u_n)}{|x|^\beta}dx\ra
   c\,\,{\rm as}\,\,n\ra\infty,\\\label{02}
   &&\le|\int_{\mathbb{R}^N}\le(|\nabla u_n|^{N-2}\nabla u_n\nabla\psi+V|u_n|^{N-2}u_n\psi\ri)dx
   -\int_{\mathbb{R}^N}\f{f(x,u_n)}{|x|^\beta}\psi
   dx\ri|\leq \tau_n\|\psi\|_E\qquad\quad
   \eea
    for all $\psi\in E$, where $\tau_n\ra 0$ as $n\ra\infty$. By Lemma 3.4,
   up to a subsequence, there holds
   \be\label{Fcon}\le\{\begin{array}{lll}
   u_n\rightharpoonup u \,\,{\rm weakly\,\, in}\,\, E\\[1.5ex]  u_n\ra
   u \,\,{\rm strongly\,\, in}\,\, L^q(\mathbb{R}^N),\,\,\forall q\geq
   1\\[1.5ex]
   \lim\limits_{n\ra\infty}\int_{\mathbb{R}^N}
   \f{F(x,u_n)}{|x|^\beta}dx=\int_{\mathbb{R}^N}\f{F(x,u)}{|x|^\beta}dx\\
   [1.5ex] u\,\,{\rm is\,\, a\,\,weak\,\,solution\,\,of}\,\,(\ref{prob-0}).\end{array}
   \ri.\ee
   Now suppose by contradiction $u\equiv0$. Since $F(x,0)= 0$ for all $x\in \mathbb{R}^N$,  it follows from
   (\ref{01}) and (\ref{Fcon}) that
   \be\label{N}\lim_{n\ra\infty}\|u_n\|_E^N=Nc>0.\ee
   Thanks to the hypothesis $(H_5)$, we have $0<c<\f{1}{N}\le(\f{N-\beta}{N}\f{\alpha_N}{\alpha_0}\ri)^{N-1}$
   by applying Lemma 3.3.
   Thus there exists some $\eta_0>0$ and $K>0$ such that
   $\|u_n\|_E^N\leq \le(\f{N-\beta}{N}\f{\alpha_N}{\alpha_0}-\eta_0\ri)^{N-1}$ for all $n>K$.
    Choose $q>1$ sufficiently close to $1$ such
   that $q\alpha_0\|u_n\|_E^{\f{N}{N-1}}\leq (1-\beta/N)\alpha_N-\alpha_0\eta_0/2$ for all
   $n>N$. By $(H_1)$,
   $$|f(x,u_n)u_n|\leq b_1|u_n|^N+b_2|u_n|\zeta\le(N,\alpha_0|u_n|^{\f{N}{N-1}}\ri),$$
   where the function $\zeta(\cdot,\cdot)$ is defined by (\ref{MF}).
   It follows from the H\"older inequality, Lemma 2.1 and Theorem
   A that
   \bna
   \int_{\mathbb{R}^N}\f{|f(x,u_n)u_n|}{|x|^\beta}dx&\leq&b_1\int_{\mathbb{R}^N}
   \f{|u_n|^N}{|x|^\beta}dx+b_2\int_{\mathbb{R}^N}\f{|u_n|\zeta\le(N,\alpha_0|u_n|^{\f{N}{N-1}}\ri)}{|x|^\beta}dx\\
   &\leq&b_1\int_{\mathbb{R}^N}\f{|u_n|^N}{|x|^\beta}dx+b_2\le(\int_{\mathbb{R}^N}\f{|u_n|^{q^\prime}}{|x|^\beta}dx\ri)^{1/{q^\prime}}
   \le(\int_{\mathbb{R}^N}\f{\zeta\le(N,q\alpha_0|u_n|^{\f{N}{N-1}}\ri)}{|x|^\beta}dx\ri)^{1/{q}}\\
   &\leq&b_1\int_{\mathbb{R}^N}\f{|u_n|^N}{|x|^\beta}dx+C\le(\int_{\mathbb{R}^N}\f{|u_n|^{q^\prime}}{|x|^\beta}dx\ri)^{1/{q^\prime}}
   \ra 0\quad{\rm as}\quad n\ra\infty.
   \ena
   Here we used (\ref{Fcon}) again (precisely $u_n\ra u$ in
   $L^s(\mathbb{R}^N)$ for all $s\geq 1$) in the last step of the above estimates.
    Inserting this into (\ref{02}) with $\psi=u_n$, we have
   $$\|u_n\|_E\ra 0\quad{\rm as}\quad n\ra\infty,$$
   which contradicts (\ref{N}). Therefore $u\not\equiv 0$ and we obtain a nontrivial
   weak solution of (\ref{prob-0}). $\hfill\Box$\\

\section{Multiplicity results\\}

In this section we will prove Theorem 1.2. The proof is divided into three steps, namely\\

\noindent{\bf Step 1}. {\it Let $\epsilon_1$ be given by $(ii)$ of
Lemma 3.1, and $\epsilon^*$, $\delta^*$ be given by Lemma 3.3. Then when $0<\epsilon<\epsilon_1$, there exists a
sequence $(v_n)\subset E$ such that
\be\label{mount}J_{\beta,\epsilon}(v_n)\ra c_M,\quad
J_{\beta,\epsilon}^\prime(v_n)\ra 0,\ee where $c_M$ is a min-max
value of $J_{\beta,\epsilon}$. Let $\epsilon_2=\min\{\epsilon_1,\epsilon^*\}$. Then when
$0<\epsilon<\epsilon_2$, we can take $c_M$ such that
\be\label{c-m}0<c_M<\f{1}{N}\le(\f{N-\beta}{N}\f{\alpha_N}{\alpha_0}\ri)^{N-1}-\delta^*.\ee
In addition, up to a subsequence, there holds $v_n\rightharpoonup
u_M$ weakly in $E$, and $u_M$ is a
weak solution of (\ref{prob-1}).}\\

\noindent{\it Proof}. By $(i)$ and $(ii)$ of Lemma 3.1, when
$0<\epsilon<\epsilon_1$,
 $J_{\beta,\epsilon}$ satisfies the following condition:
  $J_{\beta,\epsilon}\in
   \mathcal{C}^1(E,\mathbb{R})$; $J_{\beta,\epsilon}(0)=0$; $J_{\beta,\epsilon}(u)\geq \vartheta_\epsilon>0$ when
  $\|u\|_E=r_\epsilon$; $J_{\beta,\epsilon}(e)<0$ for some  $e\in E$  with
   $\|e\|>\max\{r_\epsilon,1\}$.
    Then using the mountain-pass theorem
   without the Palais-Smale condition \cite{Rabin}, we can find a sequence
   $(v_n)$ in $E$ such that
   $$J_{\beta,\,\epsilon}(v_n)\ra c_M>0,\quad J^\prime_{\beta,\,\epsilon}(v_n)\ra 0\,\,{\rm in}\,\, E^*,$$
   where
   $$c_M=\min_{\gamma\in\Gamma}\max_{u\in\gamma}J_{\beta,\epsilon}(u)\geq \vartheta_\epsilon$$
   is a min-max value of $J_{\beta,\,\epsilon}$, where $\Gamma=\{\gamma\in\mathcal{C}([0,1],E): \gamma(0)=0,
   \gamma(1)=e\}$. Clearly (\ref{c-m}) follows from Lemma 3.3. The last assertion follows from  Lemma 3.4 immediately.
   $\hfill\Box$\\

  \noindent{\bf Step 2}. {\it Let $r_\epsilon$ be given by $(ii)$ of Lemma 3.1 such that
  $r_\epsilon\ra 0$ as $\epsilon\ra 0$.
  There exists $\epsilon_3: 0<\epsilon_3<\epsilon_2$ such that if $0<\epsilon<\epsilon_3$,
  then there exists a sequence $(u_n)\subset E$ such that
  \be\label{neg1}J_{\beta,\epsilon}(u_n)\ra c_\epsilon:=\inf_{\|u\|_E\leq r_\epsilon}J_{\beta,\epsilon}(u)\ee
  and
  \be\label{neg2}J_{\beta,\epsilon}^\prime(u_n)\ra 0\quad{\rm in}\quad E^*\quad{\rm as}\quad n\ra\infty,\ee
  where  $c_\epsilon<0$ and $c_\epsilon\ra 0$ as $\epsilon\ra 0$.
  In addition, up to a subsequence, there holds $u_n\ra u_0$ strongly in $E$, and $u_0$ is a
weak solution of (\ref{prob-1}) with $J_{\beta,\epsilon}(u_0)=c_\epsilon$.} \\

\noindent{\it Proof}. Let $r_\epsilon$ be given by $(ii)$ of Lemma
3.1, i.e. $J_{\beta,\epsilon}(u)>\vartheta_\epsilon>0$ for all $u$
with $\|u\|_E=r_\epsilon$. Since $r_\epsilon\ra 0$ as $\epsilon\ra
0$, one can choose $\epsilon_3: 0<\epsilon_3<\epsilon_2$ such that
when $0<\epsilon<\epsilon_3$,
\be\label{r-e}r_\epsilon<\le(\f{N-\beta}{N}\f{\alpha_N}{\alpha_0}\ri)^{\f{N-1}{N}}.\ee
By $(H_1)$ and $(H_2)$, we have \be\label{F}F(x,u)\leq b_1
|u|^N+b_2|u|\zeta
\le(N,\alpha_0\|u\|_E^{{N}/{(N-1)}}\le({|u|}/{\|u\|_E}\ri)^{N/(N-1)}\ri).\ee
Here again $\zeta(\cdot,\cdot)$ is defined by (\ref{MF}). When
$\|u\|_E\leq r_\epsilon$, we have
$\alpha_0\|u\|_E^{N/(N-1)}<(1-\beta/N)\alpha_N$. It then follows
from Lemma 2.1 and Theorem A that $F(x,u)/|x|^\beta$ is bounded in
$L^p(\mathbb{R}^N)\cap L^1(\mathbb{R}^N)$ for some $p>1$ when $\|u\|_E\leq r_\epsilon$.
Hence $J_{\beta,\epsilon}$ has lower bound on the ball $B_{r_\epsilon}=\{u\in E:
\|u\|_E\leq r_\epsilon\}$.

Since the closure of $B_{r_\epsilon}$,
$\overline{B}_{r_\epsilon}\subset E$ is a complete metric space with
the metric given by the norm of $E$, convex and $J_{\beta,\epsilon}$
is of class $\mathcal{C}^1$ and has lower bound on
$\overline{B}_{r_\epsilon}$. By the Ekeland's variational principle
\cite{Mawin}, there exists a sequence $(u_n)\subset
\overline{B}_{r_\epsilon}$ such that (\ref{neg1}) and (\ref{neg2})
hold.

By $(iii)$ of Lemma 3.1, $c_\epsilon<0$. Since $r_\epsilon\ra 0$ as
$\epsilon\ra 0$, noticing (\ref{F}), we have by using the H\"older
inequality and Lemma 2.4
$$\sup_{\|u\|_E\leq r_\epsilon}\int_{\mathbb{R}^N}\f{F(x,u)}{|x|^\beta}dx\ra 0,\quad
\sup_{\|u\|_E\leq r_\epsilon}\int_{\mathbb{R}^N}h udx\ra 0 $$ as
$\epsilon\ra 0$. This implies $c_\epsilon\ra 0$ as $\epsilon\ra 0$.

Now we are proving the last assertion. Assume $u_n\rightharpoonup
u_0$ weakly in $E$. (\ref{neg2}) is equivalent to
\be\label{eq}|\langle
  J_{\beta,\epsilon}^\prime(u_n),\phi\rangle|\leq \tau_n\|\phi\|_E,\quad\forall \phi\in E,\ee
  where $\tau_n\ra 0$ as $n\ra\infty$. Recalling (\ref{j1}) and choosing $\phi=u_n-u_0$ in (\ref{eq}), we have
  \bna\int_{\mathbb{R}^N}\le(|\nabla u_n|^{N-2}\nabla u_n\nabla(u_n-u_0)+V(x)|u_n|^{N-2}u_n(u_n-u_0)
  \ri)dx\\-\int_{\mathbb{R}^N}\f{f(x,u_n)}{|x|^\beta}(u_n-u_0)dx-\epsilon\int_{\mathbb{R}^N}h(u_n-u_0)dx=o_n(1),\ena
  where $o_n(1)\ra 0$ as $n\ra\infty$. H\"older inequality together with (\ref{r-e}), Theorem A and Lemma 2.4 implies that
  $$\int_{\mathbb{R}^N}\f{f(x,u_n)}{|x|^\beta}(u_n-u_0)dx=o_n(1),\quad \epsilon\int_{\mathbb{R}^N}h(u_n-u_0)dx=o_n(1).$$
  Hence
  \be\label{J-e}\int_{\mathbb{R}^N}\le(|\nabla u_n|^{N-2}\nabla u_n\nabla(u_n-u_0)+V(x)|u_n|^{N-2}u_n(u_n-u_0)
  \ri)dx=o_n(1).\ee
  On the other hand, since $u_n\rightharpoonup u_0$ weakly in $E$, we obtain
  \be\label{J-2}\int_{\mathbb{R}^N}\le(|\nabla u_0|^{N-2}\nabla u_0\nabla(u_n-u_0)+V(x)|u_0|^{N-2}u_0(u_n-u_0)
  \ri)dx=o_n(1).\ee
  Subtracting (\ref{J-2}) from (\ref{J-e}), using a well known inequality (see for example Chapter 10 of \cite{Lind})
 \be\label{ineq}2^{2-N}| b- a|^N\leq \langle|b|^{N-2}b-|a|^{N-2}a, b- a\rangle,\quad\forall a,\,b\in
 \mathbb{R}^N,\ee
 we obtain $\|u_n-u_0\|_E^N\ra 0$ and thus $u_n\ra u_0$ strongly in $E$ as $n\ra\infty$. Since $J_{\beta,\epsilon}\in\mathcal{C}^1(E,\mathbb{R})$,
 there hold $J_{\beta,\epsilon}(u_0)=c_\epsilon$ and $J_{\beta,\epsilon}^\prime(u_0)=0$, i.e. $u_0$ is a weak solution of
 (\ref{prob-1}). $\hfill\Box$\\

 \noindent{\bf Step 3}. {\it There exists $\epsilon_0: 0<\epsilon_0<\epsilon_3$ such that if $0<\epsilon<\epsilon_0$,
 then $u_M\not\equiv u_0$.} \\

 \noindent{\it Proof}.  Suppose by contradiction that $u_M\equiv u_0$. Then $v_n\rightharpoonup u_0$ weakly in $E$. By (\ref{mount}),
 \be\label{mt} J_{\beta,\epsilon}(v_n)\ra c_M>0,\quad
  |\langle  J_{\beta,\epsilon}^\prime(v_n),\phi\rangle|\leq \gamma_n\|\phi\|_E\ee
  with $\gamma_n\ra 0$ as $n\ra\infty$. On one hand, by Lemma 3.4, we have
 \be\label{cv}\int_{\mathbb{R}^N}\f{F(x,v_n)}{|x|^\beta}dx\ra\int_{\mathbb{R}^N}\f{F(x,u_0)}{|x|^\beta}dx
 \quad{\rm as}\quad n\ra\infty.\ee
 Here and in the sequel, we do not distinguish sequence and subsequence. On the other hand,
 since $v_n\rightharpoonup u_0$ weakly in $E$, it follows from
 the H\"older inequality and Lemma 2.4 that
 \be\label{ch}
 \int_{\mathbb{R}^N}hv_ndx\ra \int_{\mathbb{R}^N}hu_0dx.\quad{\rm as}\quad n\ra\infty.
 \ee
 Inserting (\ref{cv}) and (\ref{ch}) into (\ref{mt}), we obtain
 \be\label{cM}
 \f{1}{N}\|v_n\|_E^N=c_M+\int_{\mathbb{R}^N}\f{F(x,u_0)}{|x|^\beta}dx+\epsilon\int_{\mathbb{R}^N}hu_0dx+o_n(1).
 \ee
 In the same way, one can derive
 \be\label{c0}
 \f{1}{N}\|u_n\|_E^N=c_\epsilon+\int_{\mathbb{R}^N}\f{F(x,u_0)}{|x|^\beta}dx+\epsilon\int_{\mathbb{R}^N}hu_0dx+o_n(1).
 \ee
 Combining (\ref{cM}) and (\ref{c0}), we have
 \be\label{Lp}
 \|v_n\|_E^N-\|u_0\|_E^N=N\le(c_M-c_\epsilon+o_n(1)\ri).
 \ee
 From Step 2, we know that $c_\epsilon\ra 0$ as $\epsilon\ra 0$. This together with (\ref{c-m}) leads to the existence
 of $\epsilon_0: 0<\epsilon_0<\epsilon_3$ such that if $0<\epsilon<\epsilon_0$, then
 \be\label{diff}0<c_M-c_\epsilon<\f{1}{N}\le(\f{N-\beta}{N}\f{\alpha_N}{\alpha_0}\ri)^{N-1}.\ee
  Write
 $$w_n=\f{v_n}{\|v_n\|_E},\quad w_0=\f{u_0}{\le(\|u_0\|_E^N+N(c_M-c_\epsilon)\ri)^{1/N}}.$$
 It follows from (\ref{Lp}) and $v_n\rightharpoonup u_0$ weakly in
 $E$ that
 $w_n\rightharpoonup w_0$ weakly in $E$. Notice that
 $$\int_{\mathbb{R}^N}\f{\zeta\le(N,\alpha_0|v_n|^{N/(N-1)}\ri)}{|x|^\beta}dx=\int_{\mathbb{R}^N}
 \f{\zeta\le(N,\alpha_0\|v_n\|_E^{{N}/{(N-1)}}|w_n|^{N/(N-1)}\ri)}{|x|^\beta}dx.$$
 By (\ref{Lp}) and (\ref{diff}), a straightforward calculation shows
 $$\lim_{n\ra\infty}\alpha_0\|v_n\|_E^{\f{N}{N-1}}\le(1-\|w_0\|_E^N\ri)^{\f{1}{N-1}}<
 \le(1-\f{\beta}{N}\ri)\alpha_N.$$
 Whence Lemma 2.3 together with Lemma 2.1 implies that
 ${\zeta\le(N,\alpha_0|v_n|^{N/(N-1)}\ri)}/{|x|^\beta}$
  is bounded in $L^q(\mathbb{R}^N)$ for some
 $q>1$. By $(H_1)$,
 $$|f(x,v_n)|\leq b_1|v_n|^{N-1}+b_2\zeta(N,\alpha_0|v_n|^{\f{N}{N-1}}).$$
 Then it follows from the continuous embedding $E\hookrightarrow L^p(\mathbb{R}^N)$
 for all $p\geq 1$ that
 $f(x,v_n)/|x|^\beta$ is bounded in $L^{q_1}(\mathbb{R}^N)$ for
 some $q_1$: $1<q_1<q$.
  This together with Lemma 2.4 and the H\"older inequality gives
 \be\label{f-c}
 \le|\int_{\mathbb{R}^N}\f{f(x,v_n)(v_n-u_0)}{|x|^\beta}dx\ri|\leq
 \le\|\f{f(x,v_n)}{|x|^\beta}\ri\|_{L^{q_1}(\mathbb{R}^N)}\le\|v_n-u_0\ri\|_{L^{{q_1^\prime}}(\mathbb{R}^N)}
 \ra 0,
 \ee
 where $1/q_1+1/q_1^\prime=1$.

 Taking $\phi=v_n-u_0$ in (\ref{mt}), we have by using (\ref{ch}) and (\ref{f-c}) that
 \be\label{cc}
  \int_{\mathbb{R}^N}\le(|\nabla v_n|^{N-2}\nabla v_n\nabla(v_n-u_0)+V(x)|v_n|^{N-2}v_n(v_n-u_0)\ri)dx
  \ra 0.
 \ee
 However the fact $v_n\rightharpoonup u_0$ weakly in $E$ leads to
 \be\label{cc1}
  \int_{\mathbb{R}^N}\le(|\nabla u_0|^{N-2}\nabla u_0\nabla(v_n-u_0)+
  V(x)|u_0|^{N-2}u_0(v_n-u_0)\ri)dx
  \ra 0.
 \ee
 Subtracting (\ref{cc1}) from (\ref{cc}), using the inequality (\ref{ineq}), we have
   $$\|v_n-u_0\|_E^N\ra
 0.$$ 
  This together with (\ref{Lp}) implies that 
 $$c_M=c_\epsilon,$$
 which is absurd since $c_M>0$ and $c_\epsilon<0$. Therefore we end Step 3 and complete the proof of Theorem 1.2.
 $\hfill\Box$\\

 {\bf Acknowledgement.} The author is supported by the program for
 NCET.

\end{document}